\documentclass{agtart_a}
\pdfoutput=1

\usepackage{pinlabel}

\title{Non-isotopic Heegaard splittings of Seifert fibered spaces}

\author{David Bachman}
\givenname{David}
\surname{Bachman}
\address{Mathematics Department\\
Pitzer College\\\newline
1050 North Mills Avenue\\
Claremont CA 91711\\USA}
\email{bachman@pitzer.edu}

\author{Ryan Derby-Talbot}
\givenname{Ryan}
\surname{Derby-Talbot}
\address{Mathematics Department\\
The University of Texas at Austin\\\newline
Austin TX 78712-0257\\
USA}
\email{rdtalbot@math.utexas.edu}

\def\appendixauthor#1{\author{Appendix by #1}}
\appendixauthor{Richard Weidmann}
\givenname{Richard}
\surname{Weidmann}
\address{{\rm R Weidmann:} Fachbereich Informatik und Mathematik\\\newline
Johann Wolfgang Goethe-Universit\"at\\
60054 Frankfurt\\Germany}
\email{rweidman@math.uni-frankfurt.de}



\received{2 May 2005}
\revised{6 December 2005}
\accepted{6 February}
\proposed{}
\seconded{}
\publishedonline{12 March 2006}
\published{12 March 2006}

\volumenumber{6}
\issuenumber{}
\publicationyear{2006}
\papernumber{12}
\startpage{351}
\endpage{372}

\doi{}
\MR{}
\Zbl{}

\arxivreference{math.GT/0504605}

\AtBeginDocument{\let\bar\wbar\let\tilde\wtilde\let\hat\what%
}
\swapnumbers

\makeatletter
\def\cnewtheorem#1[#2]#3{\newtheorem{#1}{#3}[section]
\expandafter\let\csname c@#1\endcsname\c@pro}

\cnewtheorem{thm}[pro]{Theorem}
\cnewtheorem{lem}[pro]{Lemma}
\cnewtheorem{clm}[pro]{Claim}
\cnewtheorem{example}[pro]{Example}
\cnewtheorem{cnj}[pro]{Conjecture}
\cnewtheorem{cor}[pro]{Corollary}
\cnewtheorem{quest}[pro]{Question}

\theoremstyle{definition}
\cnewtheorem{dfn}[pro]{Definition}

\theoremstyle{remark}
\newtheorem*{rmk}{Remark}

\makeatother  

\newcommand{\VV}{\mathcal V}
\newcommand{\WW}{\mathcal W}
\newcommand{\Hts}{H_{(t,s)}}

\begin{document}

\begin{asciiabstract}
We find a geometric invariant of isotopy classes of strongly
irreducible Heegaard splittings of toroidal 3-manifolds. Combining
this invariant with a theorem of R Weidmann, proved here in the
appendix, we show that a closed, totally orientable Seifert fibered
space M has infinitely many isotopy classes of Heegaard splittings
of the same genus if and only if M has an irreducible, horizontal
Heegaard splitting, has a base orbifold of positive genus, and is not
a circle bundle. This characterizes precisely which Seifert fibered
spaces satisfy the converse of Waldhausen's conjecture.
\end{asciiabstract}

\begin{abstract}
We find a geometric invariant of isotopy classes of strongly
irreducible Heegaard splittings of toroidal 3--manifolds. Combining
this invariant with a theorem of R Weidmann, proved here in the
appendix, we show that a closed, totally orientable Seifert fibered
space $M$ has infinitely many isotopy classes of Heegaard splittings
of the same genus if and only if $M$ has an irreducible, horizontal
Heegaard splitting, has a base orbifold of positive genus, and is not
a circle bundle. This characterizes precisely which Seifert fibered
spaces satisfy the converse of Waldhausen's conjecture.
\end{abstract}

\maketitle

\section{Introduction}

The recent proof of Waldhausen's conjecture (Li \cite{li}) (see also
work of Jaco and Rubinstein \cite{jr1,jr2}) establishes that a
3--manifold $M$ admits infinitely many non-isotopic Heegaard splittings
of some genus only if $M$ contains an incompressible torus. We are
interested in the converse of this statement. The only known examples
of 3--manifolds that admit infinitely many non-isotopic Heegaard
splittings of the same genus are given by Morimoto and Sakuma
\cite{sakuma,ms:91}. However, these examples are somewhat special. In
this paper, we give a complete characterization of closed, totally
orientable Seifert fibered spaces that satisfy the converse of
Waldhausen's conjecture.

In light of Li's result one would expect to use an essential torus when trying to distinguish isotopy classes of Heegaard splittings. To this end we have the following result, which is a weak version of \fullref{t:main}.

\medskip
\noindent {\bf \fullref{t:main}$'$}\qua {\sl Let $T$ be an
essential torus in an irreducible 3--manifold $M$. Suppose $H$ is a
strongly irreducible Heegaard surface in $M$ whose minimal essential
intersection number with $T$ is greater than two, and $H'$ is any
other Heegaard surface in $M$. If $H$ and $H'$ meet $T$ in different
slopes then they are not isotopic.}
\medskip

The term {\it essential intersection number} refers to the value of $|H \cap T|$ when the two surfaces are isotoped to meet in a collection of loops that are essential on both. It is well known that any strongly irreducible Heegaard surface can be isotoped to meet any essential surface in such a fashion.

Our primary goal is to distinguish non-isotopic splittings of Seifert fibered spaces. In this context we prove the following strengthening of \fullref{t:main}$'$:

\medskip
\noindent {\bf \fullref{t:SFSisotopy}}\qua {\sl Let $M$ be a
closed, totally orientable Seifert fibered space which is not a circle
bundle with Euler number $\pm 1$. Let $H$ be a strongly irreducible
Heegaard surface in $M$ and $T$ be a non-separating, vertical,
essential torus. Then the isotopy class of $H$ determines at most two
slopes on $T$.}
\medskip

In particular, if three strongly irreducible Heegaard surfaces in such a Seifert fibered space meet some essential torus in different slopes then at most two of them are isotopic. This result is stronger than \fullref{t:main}$'$ because there is no assumption on how many times any of these Heegaard surfaces meets the torus $T$. 

\fullref{t:SFSisotopy} leaves open the possibility that a circle bundle over a surface may admit an irreducible Heegaard splitting that can be isotoped to meet some vertical essential torus in infinitely many slopes. The appendix, by R Weidmann, includes a proof that this phenomenon does happen. Moreover, the Heegaard splitting in this case is unique: 

\begin{thm}[Weidmann]
\label{t:CircleBundles}
Suppose $M$ is an orientable circle bundle over an orientable surface of positive genus. Then $M$ admits a unique irreducible Heegaard splitting up to isotopy. 
\end{thm}

In addition to this Weidmann proves in the appendix an algebraic analog of the above theorem when the Euler number is $\pm 1$: Nielsen equivalence classes of the generating sets of the fundamental group of such a  manifold are equivalent. Interestingly, the algebraic formulation of this theorem motivates his topological argument used to establish \fullref{t:CircleBundles}.

The following characterization of Seifert fibered spaces that contain an infinite collection of non-isotopic Heegaard splittings of some genus now follows from Theorems \ref{t:SFSisotopy} and \ref{t:CircleBundles}.

\begin{thm}
\label{t:SFSclassification}
Let $M$ be a closed, totally orientable Seifert fibered space. Then $M$ admits infinitely many non-isotopic Heegaard splittings of some genus if and only if 
 \begin{enumerate}
 	\item $M$ has at least one irreducible, horizontal Heegaard splitting,
	\item $M$ has a base orbifold with positive genus, and
	\item $M$ is not a circle bundle. 
\end{enumerate}
\end{thm}

See \fullref{s:Definitions} below for the relevant definitions. 

\begin{proof}
Moriah and Schultens have shown that irreducible Heegaard splittings
of totally orientable Seifert fibered spaces are either vertical or
horizontal \cite{mr:98}. It follows from this classification and
results of Lustig and Moriah \cite{lm:91} and Schultens
\cite{schultens:96} that a Seifert fibered space can admit infinitely
many non-isotopic Heegaard splittings of some genus only if it admits
an irreducible, horizontal Heegaard splitting. Precisely which Seifert
fibered spaces have irreducible horizontal Heegaard splittings have
been classified by Sedgwick in terms of the Seifert data
\cite{sedgwick:99}. In particular, however, note that Moriah and
Schultens had previously shown that a circle bundle can only admit an
irreducible horizontal Heegaard splitting if its Euler number is $\pm
1$ (see  \cite[Corollary 0.5]{mr:98}).

An understanding of horizontal Heegaard splittings reveals that any infinite collection must be obtained by Dehn twists in vertical tori (see, for example, Hatcher's proof that incompressible surfaces in Seifert fibered spaces are either vertical or horizontal \cite{hatcher}). So the question of whether a given closed, totally orientable Seifert fibered space admits an infinite collection of non-isotopic splittings of some genus is reduced to determining when Dehn twisting a horizontal splitting in a vertical torus produces a non-isotopic splitting. This is recognized by Sedgwick in the following:

\begin{quote}
The author suspects ... that some Seifert fibered spaces will posses
 an infinite number of non-isotopic but homeomorphic splittings
 obtained by twisting a given horizontal splitting in vertical tori
 \cite[page 178, line -2]{sedgwick:99}.  \end{quote}

Suppose, then, that $\mathcal V \cup _H \mathcal W$ is an irreducible,
horizontal Heegaard splitting of a Seifert fibered space $M$ and $T$
is a vertical torus. Assume first that $T$ separates $M$ into $X$ and
$Y$. Then $T$ separates $H$ into a horizontal surface $H_X$ in $X$
(say) and a surface which is not horizontal in $Y$. But then $H_X$ is
a union of fibers in a fibration of $X$ over $S^1$ (see Jaco
\cite[Theorem VI.34]{jaco:80}). Hence, the effect of Dehn twisting $H$
about $T$ can be undone by pushing $H_X$ around the fibration. The
conclusion is that a Dehn twist about a separating, vertical torus
produces an isotopic Heegaard splitting. In particular, if the base
orbifold of $M$ is a sphere then every vertical torus separates, and
hence $M$ has finitely many non-isotopic Heegaard splittings in each
genus.

Now assume the base orbifold of $M$ has positive genus. If $M$ is a circle bundle, then by \fullref{t:CircleBundles} $M$ admits finitely many Heegaard splittings, up to isotopy. Henceforth, assume $M$ is not a circle bundle. 

It follows from Theorems 2.6 and 5.1 of Moriah--Schultens \cite{mr:98}
that in Seifert fibered spaces with positive genus base orbifold, all
irreducible, horizontal Heegaard splittings are strongly
irreducible. Hence, the surface $H$ is strongly irreducible. As the
base orbifold has positive genus, we may find a pair of non-separating
vertical tori $T_1$ and $T_2$ which meet in a single fiber $f$. A
horizontal Heegaard surface such as $H$ meets each of these tori in
loops that are transverse to $f$. Dehn twisting $H$ about $T_2$ has
the same effect, on $T_1$, as Dehn twisting $H \cap T_1$ about
$f$. Hence the new splitting surface meets $T_1$ in a different slope
than the original splitting surface. Iterating the Dehn twist about
$T_2$ thus produces an infinite collection of Heegaard splittings, all
of which meet $T_1$ in distinct slopes. It now follows from 
\fullref{t:SFSisotopy} that this collection contains infinitely many
non-isotopic splittings.
\end{proof}

The authors would like to thank Cameron Gordon and Yo'av Rieck for helpful comments, and especially Richard Weidmann for providing the appendix. 

\section{Definitions}
\label{s:Definitions}

\subsection{Essential curves, surfaces and intersections}

A sphere in a 3--manifold is {\it essential} if it does not bound a ball. If a 3--manifold does not contain any essential spheres then it is said to be {\it irreducible}.

A loop $\gamma$ on a surface $F$ if called {\it inessential} if it bounds a disk in $F$ and {\it essential} otherwise. The intersection between surfaces $H$ and $T$ in a 3--manifold is {\it compression free} if the surfaces are transverse and every loop contained in their intersection is either essential or inessential on both surfaces. Their intersection  is {\it essential} if every loop contained in their intersection is essential on both.

If $T$ is a torus then a {\it slope} on $T$ is the isotopy class of an essential loop. If $H$ is some other surface then the {\it slope} of $H \cap T$ is the slope of any component of $H \cap T$ which is essential on $T$. Note that this is only defined when there is such a component of $H \cap T$.

Suppose $F$ is embedded in a 3--manifold $M$. A {\it compressing disk} for $F$ is a disk $D$ such that $D \cap F=\partial D$ is essential on $F$. A surface is {\it compressible} if there is a compressing disk for it, and {\it incompressible} otherwise. A surface of positive genus in a 3--manifold is said to be {\it essential} if it is incompressible and non-boundary parallel. 

\subsection{Heegaard splittings}

A {\it handlebody} is a 3--manifold which is homeomorphic to the neighborhood of a connected graph in $\mathbb R^3$. An expression of a 3--manifold $M$ as $\mathcal V \cup _H \mathcal W$ is called a {\it Heegaard splitting} if $\mathcal V$ and $\mathcal W$ are handlebodies. The surface $H$ is called the {\it Heegaard surface}. 

A Heegaard splitting $\mathcal V \cup _H \mathcal W$ is said to be {\it reducible} if there are compressing disks $V \subset \mathcal V$ and $W \subset \mathcal W$ for the surface $H$ such that $\partial V=\partial W$, and {\it irreducible} otherwise. A Heegaard splitting is said to be {\it weakly reducible} if there are similar disks $V$ and $W$ such that $V \cap W =\emptyset$, and {\it strongly irreducible} otherwise. 

\subsection{Seifert fibered spaces}
\label{s:SFS}

A 3--manifold $M$ is a {\it Seifert fibered space} if there is a projection map $p:M \to \mathcal O$, where $\mathcal O$ is a surface and $p^{-1}(x)$ is a circle for each $x \in \mathcal O$. The surface $\mathcal O$ is called the {\it base surface} of the fibration, and inherits from $p$ a natural structure as an orbifold. If $x$ is a cone point of $\mathcal O$ then we say $p^{-1}(x)$ is an {\it exceptional fiber}. For all other $x$ we say $p^{-1}(x)$  is a {\it regular fiber}. A Seifert fibered space $M$ is {\it totally orientable} if it is orientable and its base orbifold $\mathcal O$ is orientable.

A surface in a Seifert fibered space is {\it horizontal} if it is transverse to each fiber. The following facts are known about horizontal surfaces. See, for example, Jaco \cite{jaco:80}.

\begin{enumerate}
	\item If a Seifert fibered space contains an essential surface with non-zero Euler characteristic then it can be made horizontal. 
	\item Every Seifert fibered space with boundary contains a horizontal surface. 
	\item If a totally orientable Seifert fibered space $M$ contains a connected, horizontal surface $F$ then $M$ can be obtained from $F \times I$ by identifying $F \times \{0\}$ with $F \times \{1\}$ via some homeomorphism.
	\item \label{t:NoExcep} If a Seifert fibered space contains a horizontal surface which meets a regular fiber once, then it contains no exceptional fibers.
\end{enumerate}

A Heegaard splitting $\mathcal V \cup _H \mathcal W$ of a Seifert fibered space $M$ is said to be {\it horizontal} if the surface $H$ can be obtained by the following construction. Let $M(f)$ denote the Seifert fibered space obtained from $M$ by removing a neighborhood of some fiber $f$. Then $M(f)$ has boundary, and can therefore be obtained from some surface $F$ with connected boundary by forming $F \times I$ and identifying $F \times \{0\}$ with $F \times \{1\}$ via some homeomorphism. Now take two parallel copies of $F$ and join them by a subannulus of $\partial M(f)$ to form $H$. Let $D$ denote a meridional disk for the solid torus attached to $M(f)$ to form $M$. The surface $H$ obtained by the above procedure will be a Heegaard surface in $M$ when $\partial D$ meets $\partial F$ in a point.

\section{Sweepouts}

Let $H$ denote a Heegaard surface in a 3--manifold $M$. Then there is a {\it sweepout} of $M$ by surfaces parallel to $H$. To be precise, there is a pair of graphs $\Sigma _0$ and $\Sigma _1$ embedded in $M$ and a continuous map $\Phi :H \times I \to M$ such that
	\begin{itemize}
		\item $\Phi (H  \times \{0\})=\Sigma _0$,
		\item $\Phi (H \times \{1\})=\Sigma _1$, 
		\item there is an $s$ such that $\Phi (H \times \{s\})=H$, and
		\item $\Phi$ is a homeomorphism when restricted to $H \times (0,1)$.
	\end{itemize}
Henceforth, we denote $\Phi (H  \times \{s\})$ as $H_s$. 

Now suppose $M$ is irreducible, $T$ is an essential torus in $M$ and $H$ is strongly irreducible. The sweepout $\Phi$ induces a height function $h:T \to I$ as follows: if $x \in T \cap H_s$ then $h(x)=s$. We assume $\Phi$ is chosen so that $h$ is Morse on $h^{-1}(0,1)$. 

\begin{lem}
\label{l:EssentialIntersection}
There are values $s_- < s_+ $ corresponding to saddle tangencies such that $H_s \cap T$ is compression free if and only if $H_s$ is transverse to $T$ and $s_- \le s \le s_+$.
\end{lem}

The fact that there exists a regular value $s$ such that $H_s \cap T$
is compression free is a well known result, and is established here in
Claims \ref{c:claim1} through \ref{c:claim4} of the following
proof. The real content of \fullref{l:EssentialIntersection} is that
the closure of all $s$ such that $H_s \cap T$ is compression free is a
connected interval. This is established in \fullref{c:claim5}, which
is reminiscent of Bachman--Schleimer \cite[Claim 6.7]{fiber}.

\proof
For each $s\in (0,1)$ the surface $H_s$ separates $M$ into handlebodies $\VV_s$ and $\WW_s$, where $a < b$ implies $\VV_a \subset \VV _b$. Let $s_0=0$, $\{s_i\}_{i=1}^{n-1}$ the values of $s$ where $H_s$ is not transverse to $T$, and $s_n=1$.

We now label the intervals $[s_i,s_{i+1}]$  as follows. If, for some value of $s \in (s_i,s_{i+1})$, the intersection set $H_s \cap T$ contains a loop which is essential on $H_s$ and bounds a disk in $\VV_s$ then we label the interval $[s_i,s_{i+1}]$ with the letter ``$V$". Similarly, if the intersection set $H_s \cap T$ contains a loop which is essential on $H_s$ and bounds a disk in $\WW_s$ then we label the interval $[s_i,s_{i+1}]$ with the letter ``$W$". Note that  $H_s \cap T$ is compression free if and only if $s$ is in an unlabeled interval.

\begin{clm}
\label{c:claim1}
For every $s \in (0,1)$ the intersection $H_s \cap T$ contains a loop which is essential on $H_s$.
\end{clm}

\begin{proof}
Suppose not. Then a standard innermost disk argument would show that $T$ may be isotoped to be disjoint from $H_s$, and hence lie in a handlebody. This is a contradiction, as $T$ is incompressible.
\end{proof}

\begin{clm}
The interval $[s_0,s_1]$ is labeled ``$V$" and the interval $[s_{n-1},s_n]$ is labeled ``$W$".
\end{clm}

\begin{proof}
Choose some $s$ just larger than $s_0=0$. Then $H_s$ meets $T$ in a collection of loops which all bound disks in $\VV _s$. By the previous claim at least one of these loops is essential on $H_s$, so the interval $[s_0,s_1]$ is labeled ``$V$". A symmetric argument completes the proof.
\end{proof}

\begin{clm}
No interval is labeled with both a ``$V$" and a ``$W$".
\end{clm}

\begin{proof}
Suppose this is the case for the interval $[s_i,s_{i+1}]$. Choose some $s \in (s_i,s_{i+1})$. Then there are loops in $H_s \cap T$ bounding disks in $\VV_s$ and $\WW_s$. This contradicts the strong irreducibility of $H_s$. 
\end{proof}

\begin{clm}
\label{c:claim4}
Intervals with the labels ``$V$" and ``$W$" cannot be adjacent.
\end{clm}

\begin{proof}
Suppose $[s_{i-1},s_i]$ and $[s_i,s_{i+1}]$ are adjacent intervals with different labels. Then the surface $H_{s_i}$ meets $T$ in a saddle tangency. Let $\Omega$ denote the graph $H_{s_i} \cap T$ and $N(\Omega)$ a regular neighborhood of this graph on $H_{s_i}$. 

Without loss of generality assume the label of $[s_{i-1},s_i]$ is ``$V$". For small $\epsilon$ the intersection $H_{s_i-\epsilon} \cap T$ contains a loop bounding a disk in $\VV_{s_i-\epsilon}$ (say), so there is a loop of $\partial N(\Omega)$ bounding a disk in $\VV_{s_i}$. Similarly, $H_{s_i+\epsilon} \cap T$ contains a loop bounding a disk in $\WW_{s_i+\epsilon}$, so there is a loop $\partial N(\Omega)$ bounding a disk in $\WW_{s_i}$. As these loops are either the same or are disjoint we again contradict strong irreducibility.
\end{proof}

It follows from the preceding claims that there exists an unlabeled interval. The proof of the lemma is then complete once we establish the following:

\begin{clm}
\label{c:claim5}
The union of the unlabelled intervals is connected. 
\end{clm}

\begin{proof}
Suppose $[s_i, s_{i+1}],[s_j, s_{j+1}]$ and $[s_k, s_{k+1}]$ are intervals where $i <j<k$, $[s_i, s_{i+1}]$ and $[s_k, s_{k+1}]$ are unlabeled, and $[s_j, s_{j+1}]$ has a label. Without loss of generality assume the label of $[s_j, s_{j+1}]$ is ``$V$". Then there is a disk $V \subset \VV_s$ such that $\partial V=\alpha \subset H_s \cap T$ for some $s \in (s_j,s_{j+1})$. As $T$ is incompressible an innermost disk argument can be used to show that the loop $\alpha$ bounds a disk $V' \subset T$. 

Choose $s' \in (s_i, s_{i+1})$. As $i < j$ we have $s' < s$. We claim that $V' \cap H_{s'}$ contains a loop which is essential on $H_{s'}$. If not then an innermost disk argument would show that $V'$ can be isotoped to be disjoint from $H_{s'}$. Now let $\alpha '$ denote a loop of $V' \cap H_s$ which is innermost (on $V'$) among all loops which are essential on $H_s$ (possibly $\alpha'=\alpha$). Let $V''$ denote the subdisk of $V'$ bounded by $\alpha'$. Then an innermost disk argument shows that $V''$ can be isotoped to a compressing disk for $H_s$, while still being disjoint from $H_{s'}$. As the region between $H_s$ and $H_{s'}$ is a product it follows that $V'' \subset \WW_s$. We conclude $\alpha$ is a loop of $H_s$ bounding a compressing disk in $\VV_s$ and $\alpha '$ is a loop bounding a compressing disk in $\WW_s$, contradicting the strong irreducibility of $H_s$. 

We conclude that $V' \cap H_{s'}$ contains a loop which is essential on $H_{s'}$. Let $\beta$ denote such a loop which is innermost (on $V'$). Note that as $V' \subset T$ we have $\beta \subset H_{s'} \cap T$. Since the interior of the subdisk of $T$ bounded by $\beta$ meets $H_{s'}$ in loops that are inessential on both surfaces we may remove them by an innermost disk argument. Hence, $\beta$ bounds a compressing disk for $H_{s'}$, which must lie in either $\VV_{s'}$ or $\WW_{s'}$. In either case the interval $[s_i,s_{i+1}]$ would have had a label.

If, initially, the label of $[s_j,s_{j+1}]$ was ``$W$" we would have chosen $s' \in (s_k,s_{k+1})$ and used a symmetric argument.
\end{proof}

\section{Compression free isotopies}

\begin{dfn}
Let $T$ be an essential torus in a 3--manifold $M$. An isotopy $H \times I \to M$ is {\it compression free with respect to $T$} if, for all $t \in I$ such that $H_t$ is transverse to $T$, the intersection $H_t \cap T$ is compression free.
\end{dfn}

\begin{lem}
\label{l:CompressionFreeIsotopy}
Let $H_0$ and $H_1$ denote isotopic, strongly irreducible Heegaard surfaces in an irreducible 3--manifold $M$. Let $T$ be an essential torus in $M$. Suppose $H_i \cap T$ is compression free, for $i=0,1$. Then there is an isotopy from $H_0$ to $H_1$ which is compression free with respect to $T$. Furthermore, there is such an isotopy such that the tangencies of $H_t \cap T$ which develop are either centers, saddles, or double-saddles.
\end{lem}

\begin{proof}
Let $H_t$ denote any isotopy from $H_0$ to $H_1$. We now define a two-parameter family of Heegaard surfaces. Note that for each $t$ the surface $H_t$ defines a sweepout $\Phi _t:H_t \times I \to M$. We denote $\Phi_t(H_t \times \{s\})$ as $H_{(t,s)}$. This defines a map from $H \times I \times I$ into $M$, which we can choose to be continuous in $s$ and $t$. Furthermore, there are values $s_0$ and $s_1$ such that $H_{(0,s_0)}=H_0$ and $H_{(1,s_1)}=H_1$. 

Now suppose $T$ is an essential torus in $M$. Let $\Gamma$ denote the set of points in $I \times I$ such that $H_{(t,s)}$ is not transverse to $T$. According to Cerf theory \cite{cerf:68} we may assume $\Gamma$ is a graph with vertices of valence two and four, and for each $t$ at most one vertex of $\Gamma$ is contained in $t \times I$. We say $t$ is a {\it regular value} if there is no vertex of $\Gamma$ in $t \times I$.

Let $S$ denote the closure of the set of points $(t,s) \subset I \times I$ such that $\Hts \cap T$ is compression free. We now claim that there is a path from $(0,s_0)$ to $(1,s_1)$ in $S$. Such a path defines the desired compression free isotopy from $H_0$ to $H_1$. It may pass through edges of $\Gamma$ corresponding to center or saddle tangencies, or a valence four vertex of $\Gamma$ which will correspond to two saddle tangencies.

Let $\pi:I \times I \to I$ denote projection onto the first factor. Let $p$ and $q$ denote paths in $S \subset I \times I$ (ie, embedded intervals) such that 

\begin{enumerate}
	\item $(0,s_0) \in p$,
	\item $(1,s_1) \in q$,
	\item the lengths of $\pi(p)$ and $\pi(q)$ are maximal.
\end{enumerate}

\begin{clm}
If the sum of the lengths of $\pi(p)$ and $\pi(q)$ is greater than one then there is a path in $S$ from $(0,s_0)$ to $(1,s_1)$.
\end{clm}

\begin{proof}
In this case there is an $x \in p$ and a $y \in q$ such that $\pi(x)=\pi(y)$ is a regular value of $t$. By \fullref{l:EssentialIntersection} the subinterval $r$ of $\pi(x) \times I$ connecting $x$ to $y$ is in $S$. Let $p'$ denote the subpath of $p$ connecting $(0,s_0)$ to $x$ and $q'$ the subpath of $q$ connecting $y$ to $(1,s_1)$. Then the path $p' \cup r \cup q'$ is the desired path from $(0,s_0)$ to $(1,s_1)$.
\end{proof}

\begin{clm}
The lengths of $\pi(p)$ and $\pi(q)$ are equal to one.
\end{clm}

\begin{proof}
By way of contradiction, assume the length of $\pi(p)$ is less than one. Let $(t^*,s^*)$ denote the endpoint of $p$ which is not $(0,s_0)$. 

For each $t$ there is at most one vertex of $\Gamma$ in $t \times I$. We may thus choose an $\epsilon$ small enough so that there is at most one vertex of $\Gamma$ in the rectangle $R=[t^*-\epsilon,t^*+\epsilon] \times I$. Let $t^-=t^*-\epsilon$ and $t^+=t^*+\epsilon$. We may assume that $t^-$ and $t^+$ are regular values of $t$. Finally, as $\epsilon$ is chosen to be small we may assume that there is at most one component of $\Gamma \cap R$ which is not an arc connecting $t^- \times I$ to $t^+ \times I$.

        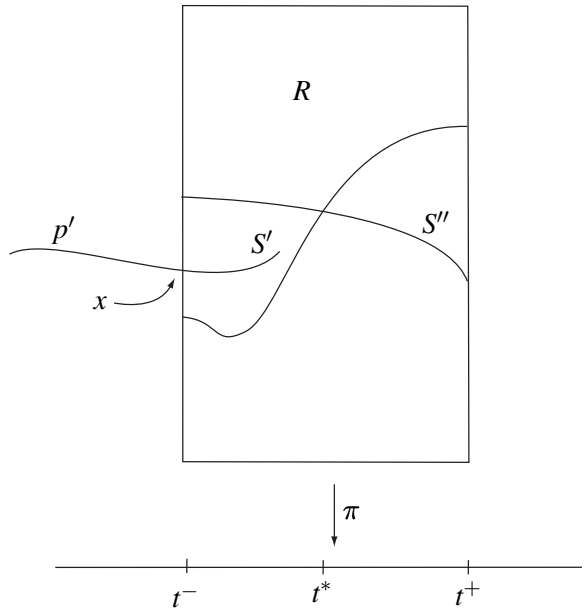
\begin{figure}[ht!]
        \labellist
        \pinlabel {$p'$}  [b]  at 117 386
        \pinlabel {$x$}   [r]  at 156 346
        \pinlabel {$R$}        at 303 513
        \pinlabel {$\pi$} [l]  at 329 181
        \pinlabel {$S'$}   <0pt,-2pt>    at 272 398
        \pinlabel {$S''$} [bl] at 393 396
        \pinlabel {$t^-$} [t] <0pt,-1pt>  at 212 133
        \pinlabel {$t^*$} [t]  at 319 133
        \pinlabel {$t^+$} [t]  at 434 133
        \pinlabel {$\pi$} [l]  at 329 181
        \endlabellist
        \begin{center}
        \includegraphics[width=3in]{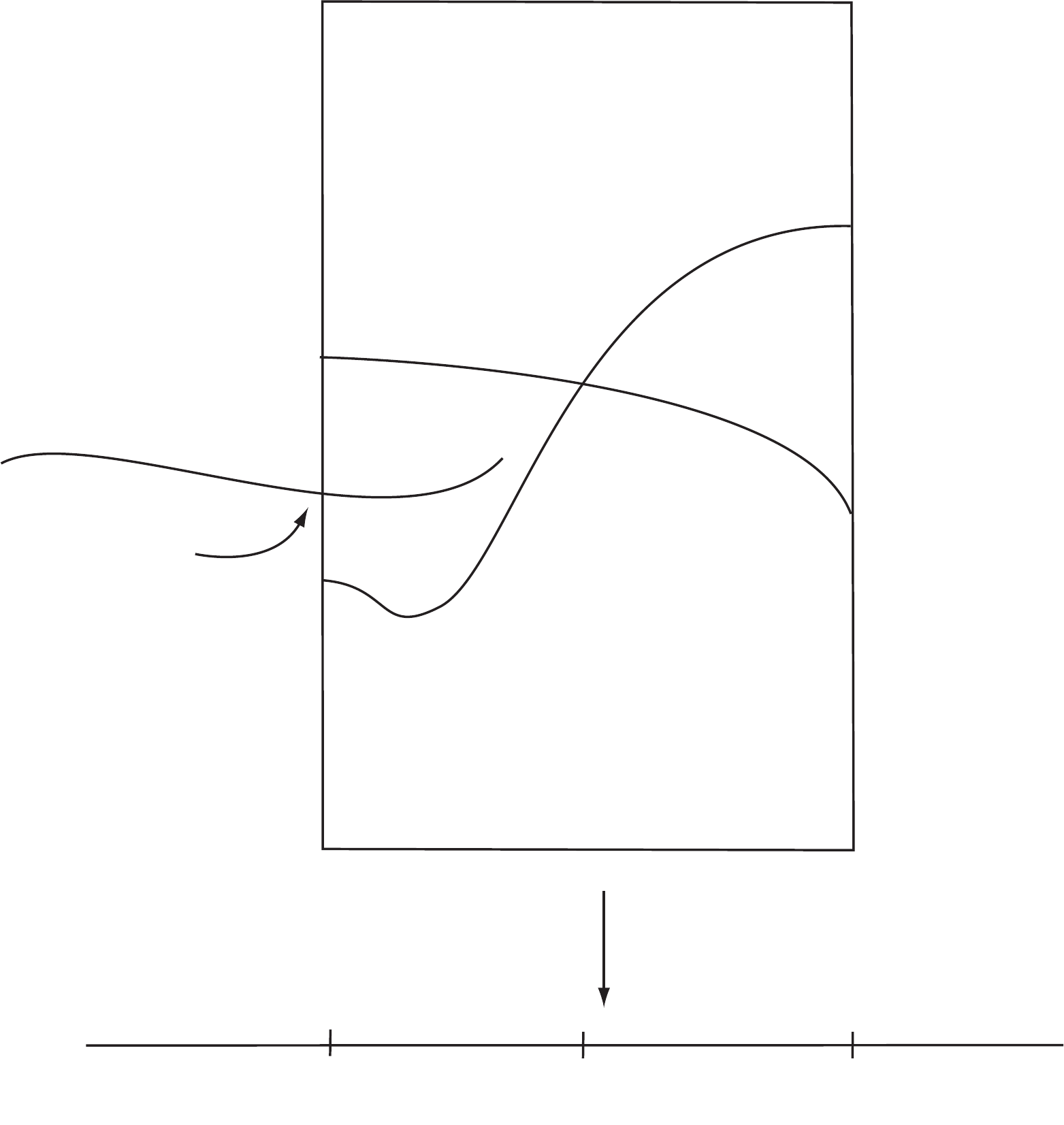}
        \caption{The rectangle $R$}
        \label{f:R}
        \end{center}
        \end{figure}

Let $p'$ denote the closure of $p \setminus R$. Let $x$ denote the endpoint of $p'$ which is not $(0,s_0)$. Note that $\pi(x)=t^-$ (see \fullref{f:R}). Let $S'$ be the closure of the component of $R \setminus \Gamma$ that contains $x$. Since $x \in S$ it follows that $S' \subset S$. If $S'$ meets the edge $t^+ \times I$ of $R$ then there is a path $p''$ in $S'$ (and hence in $S$) from $x$ to a point of $t^+ \times I$. The path $p' \cup p''$ thus contradicts the maximality of the length of $\pi(p)$. 

We assume then that $S'$ does not meet $t^+ \times I$. Let $S''$ denote the closure of a component of $R \setminus \Gamma$ which is a subset of $S$ and meets the edge $t^+ \times I$ (such a component exists by \fullref{l:EssentialIntersection}). By \fullref{l:EssentialIntersection} the set $S \cap (t^- \times I)$ is connected. Hence, if $S''$ also meets $t^- \times I$ then as before we can extend the path $p$ to $t^+ \times I$, contradicting our assumption that $\pi(p)$ is maximal.

We are now reduced to the case that $S'$ does not meet $t^+ \times I$ and $S''$ does not meet $t^- \times I$. The only way in which this can happen is if $S'$ meets $S''$ in a valence four vertex $v$ of $\Gamma$. We conclude that there is a path $p''$ which goes from $x$, through $S'$, across $v$, through $S''$, and connects to $t^+ \times I$. The path $p' \cup p''$ again contradicts the maximality of the length of $\pi(p)$. 
\end{proof}

The preceding claims complete the proof of \fullref{l:CompressionFreeIsotopy}.
\end{proof}

\begin{thm}
\label{t:main}
Let $T$ be an essential torus in an irreducible 3--manifold $M$. Suppose $H_0$ and $H_1$ are isotopic, strongly irreducible Heegaard surfaces  which meet $T$ essentially. Then either $H_0$ can be isotoped to meet a neighborhood of $T$ in a toggle or $H_0$ determines the same slope on $T$ as $H_1$. 
\end{thm}

The term {\it toggle} refers to the configuration depicted in \fullref{f:toggle}. It can be constructed as follows. Let $\alpha$ and $\beta$ be essential loops on $T$ which meet in a point $p$. Let $\Sigma$ be the graph $(\alpha \times \{0\}) \cup (p \times I) \cup (\beta \times \{1\})$ in $T \times I$. Then the frontier of a neighborhood of $\Sigma$ in $T \times I$ is a toggle. The word ``toggle" comes from the fact that such a configuration allows one to switch back and forth between two slopes in a neighborhood of $T$. 

\begin{proof}
By \fullref{l:CompressionFreeIsotopy} we know that there is a compression free isotopy from $H_0$ to $H_1$. We now discuss the various tangencies with $T$ that can develop during such an isotopy, and how they effect the slope of $H_t \cap T$. 

\medskip
\noindent {\bf Center Tangencies}\qua The simplest is a center tangency. Such tangencies only introduce or eliminate inessential loops, and hence do not change the slope of $H_t \cap T$. 

\medskip
\noindent {\bf Saddle Tangencies}\qua The next type of tangency is a saddle. If $H'$ is obtained from $H$ by passing through a saddle with $T$ then there is a disk $S$ such that $\partial S=\alpha \cup \beta$, where $S \cap T=\alpha$ and $S \cap H=\beta$ (see \fullref{f:SaddleDisk}). The surface $H'$ is then obtained from $H$ by an isotopy guided by $S$. Hence the intersection set $H' \cap T$ can be obtained from $H \cap T$ by a band sum along the arc $\alpha$. We call such a disk $S$ a {\it saddle disk}.

        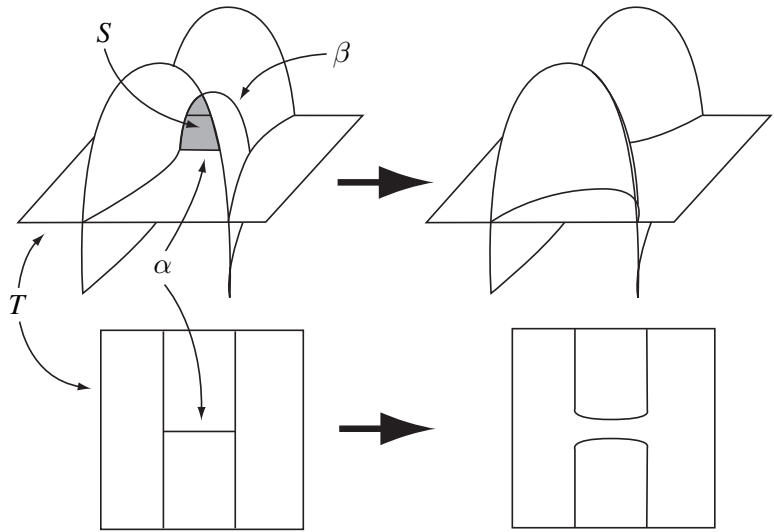
\begin{figure}[ht!]
        \labellist
        \pinlabel {$S$}  [b] at 98 541
        \pinlabel {$\alpha$} <0pt,2pt> at 134 405
        \pinlabel {$\beta$} [l] <-1pt,0pt> at 235 538
        \pinlabel {$T$}      at 45 385
        \endlabellist
        \begin{center}
        \includegraphics[width=4 in]{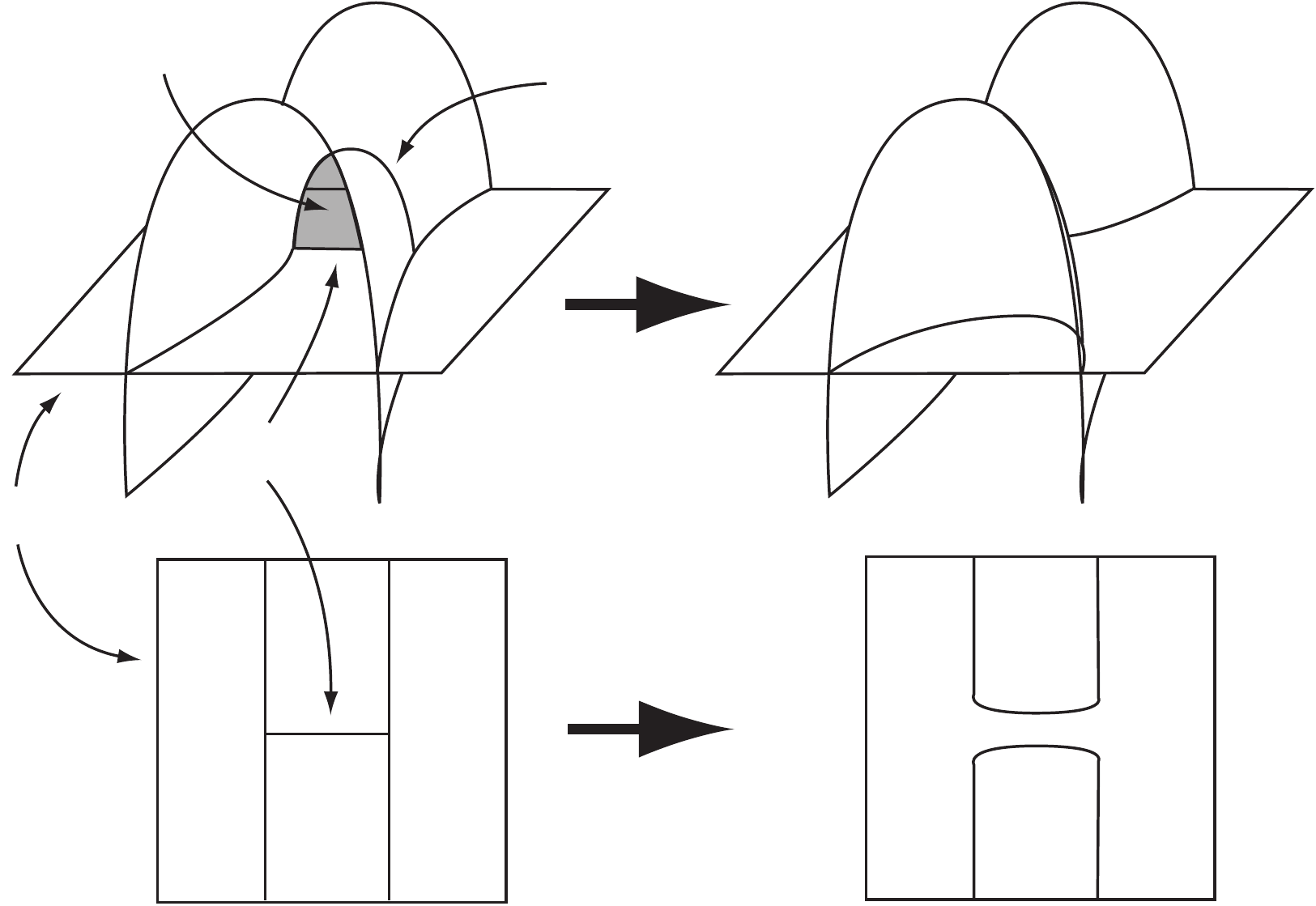}
        \caption{A saddle disk}
        \label{f:SaddleDisk}
        \end{center}
        \end{figure}

Note that the only way that the slope of $H \cap T$ can be different from the slope of $H' \cap T$ is if somehow all of the essential loops of $H \cap T$ were effected during the saddle move. But the only such loops that will be effected are those that contain the endpoints of $\alpha$. It follows that $H \cap T$ contains exactly two essential loops, and $\alpha$ is an arc which connects them. But then a band sum along $\alpha$ will produce an intersection set with no essential loops on $T$. This is impossible, as $H' \cap T$ is compression free. 

\medskip
\noindent {\bf Double-saddle Tangencies}\qua Finally we consider what happens at double-saddles. Suppose $H'$ is obtained from $H$ by passing through a double-saddle with $T$. Then there are two saddle disks $S_1$ and $S_2$, where $\partial S_i =\alpha _i \cup \beta _i$, $S_i \cap T=\alpha _i$, and $S_i \cap H=\beta _i$. The intersection set  $H' \cap T$ is obtained from $H \cap T$ by simultaneous band sums along $\alpha _1$ and $\alpha _2$.

In order for the slope of $H' \cap T$ to be different from the slope of $H \cap T$ all of the essential loops of $H \cap T$ must contain an endpoint of either $\alpha _1$ or $\alpha _2$.  This immediately implies $H \cap T$ contains at most four essential loops. The possibility that there are one or three such loops is ruled out by the fact that $H$ is separating. If there are four such loops, and each contains an endpoint of $\alpha_1$ or $\alpha _2$, then $H ' \cap T$ contains only inessential loops. This is ruled out by the fact that $H' \cap T$ is compression free. 

We conclude that if the slope of $H' \cap T$ is different from that of $H \cap T$ then $H \cap T$ contains exactly two essential loops, $\gamma _1$ and $\gamma _2$. Up to relabeling, there are now the following cases:
\begin{enumerate}
	\item $\partial \alpha _1 \subset \gamma _1$. Then a band sum along $\alpha _1$ transforms $\gamma _1$ into an essential loop $\gamma _1'$ with the same slope on $T$, and an inessential loop $\delta$. The arc $\alpha _2$ can either connect $\gamma _2$ to itself, connect $\gamma_2$ to $\delta$, or connect $\gamma _2$ to $\gamma '_1$. In the first two cases a slope change does not occur. The third case implies $H' \cap T$ contains only inessential loops, which cannot happen. 
	\item Both $\alpha_1$ and $\alpha _2$ connect $\gamma _1$ to $\gamma _2$. If $\alpha _1$ and $\alpha _2$ are on the same side of $H$ then a simultaneous band sum results in all inessential loops. We conclude $\alpha _1$ and $\alpha _2$ are are on opposite sides of $H$, as in \fullref{f:toggleT}.
\end{enumerate}

        \begin{figure}[ht!]
        \labellist
        \pinlabel {$\alpha_1$} [b] at 144 426
        \pinlabel {$\alpha_2$} [b] at 211 488
        \endlabellist
        \begin{center}
        \includegraphics[width=1 in]{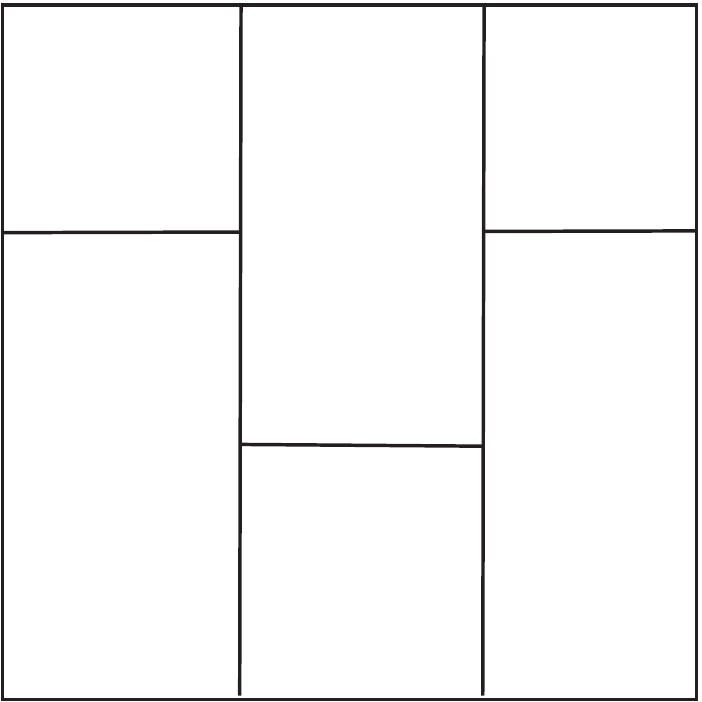}
        \caption{The set $H \cap T$ when there is a slope 
                 change at a double saddle}
        \label{f:toggleT}
        \end{center}
        \end{figure}

Now that we have narrowed down the possibilities for $H \cap T$ and $\alpha _1$ and $\alpha _2$ we must analyze the saddle disks $S_1$ and $S_2$. Before proceeding further note that if there are any inessential loops on $H \cap T$ they may be removed by an isotopy of $H$, as $H \cap T$ is compression free and $M$ is irreducible. After performing such an isotopy let $A_i$ be the annulus on $T$ bounded by $\gamma _1 \cup \gamma _2$ containing the arc $\alpha _i$. Let $D_i$ denote the disk obtained by gluing two parallel copies of $S_i$ to the disk obtained from $A_i$ by removing a neighborhood of $\alpha _i$. 

First note that if, for some $i$, the disk $D_i$ failed to be a compressing disk for $H$ it would follow that the component $H^*$ of $H \backslash T$ containing $\beta _i$ was an annulus which is parallel into $T$. Hence, a further isotopy of $H$ could push $H^*$ past $T$, removing all intersections of $H$ with $T$. As this is impossible, we conclude both $D_1$ and $D_2$ are compressing disks for $H$.

        \begin{figure}[ht!]
        \labellist
        \pinlabel {$S_1$} [t] at 228 183
        \pinlabel {$S_2$} [tr] <2pt,0pt> at 392 178
        \pinlabel {$T$} at 89 229
        \endlabellist
        \begin{center}
        \includegraphics[width=3.5 in]{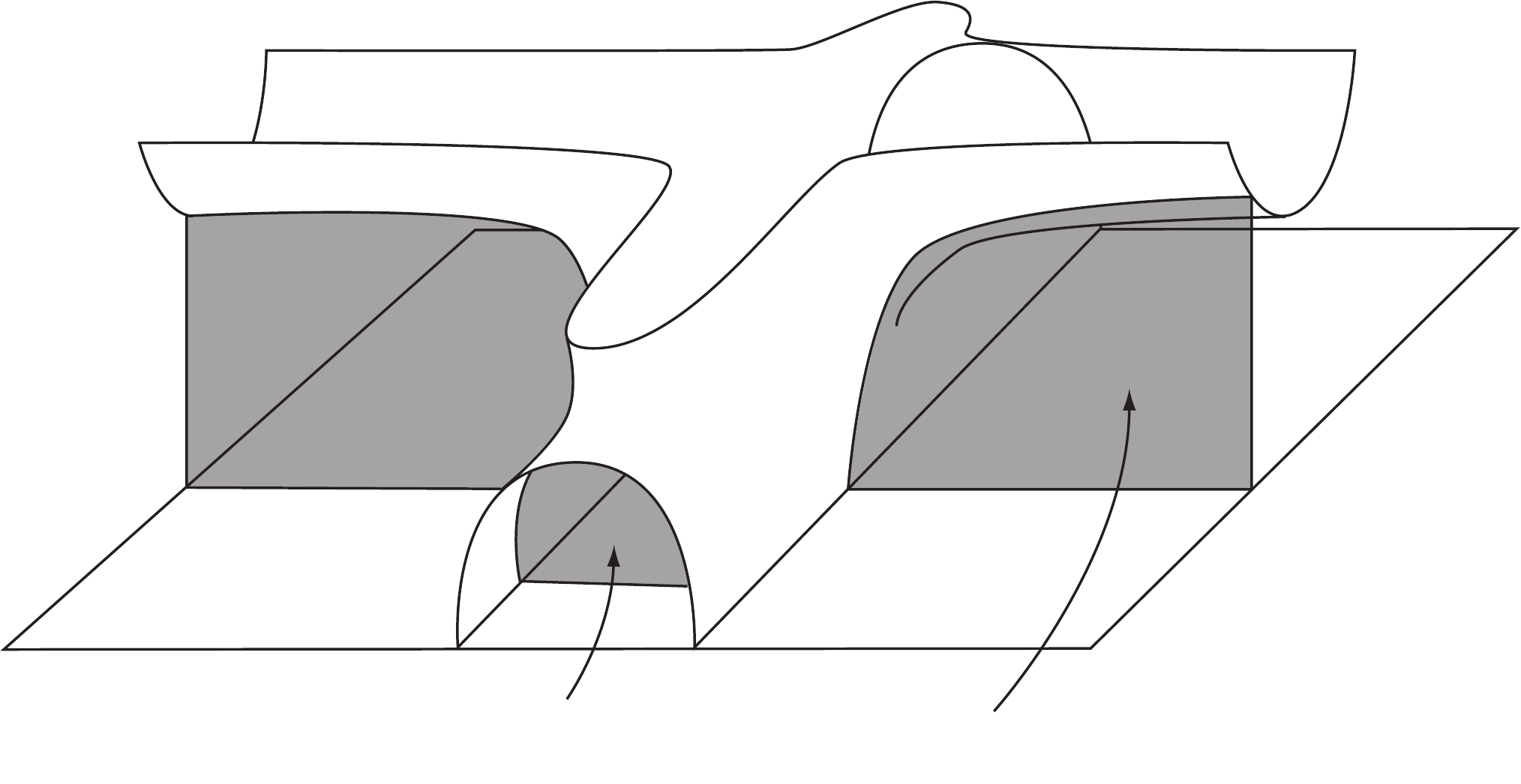}
        \caption{The surface $H \cap N(T)$}
        \label{f:S1S2}
        \end{center}
        \end{figure}

Now note that if $S_1$ and $S_2$ are incident to opposite sides of $T$ then the disks $D_1$ and $D_2$ would be disjoint. This violates the strong irreducibility of $H$. We conclude $S_1$ and $S_2$ are on the same side of $T$. Let $N(T)$ denote a copy of $T^2 \times I$ embedded in $M$ so that $T$ is the image of $T^2 \times \{0\}$. We may thus assume that $S_1$ and $S_2$ are contained in $N(T)$. This forces $H \cap N(T)$ to be as depicted in \fullref{f:S1S2}. It is now an easy exercise to see that $H \cap N(T)$ is a toggle, as depicted in \fullref{f:toggle}.

        \begin{figure}[ht!]
        \vspace{0 in}
        \begin{center}
        \includegraphics[width=2 in]{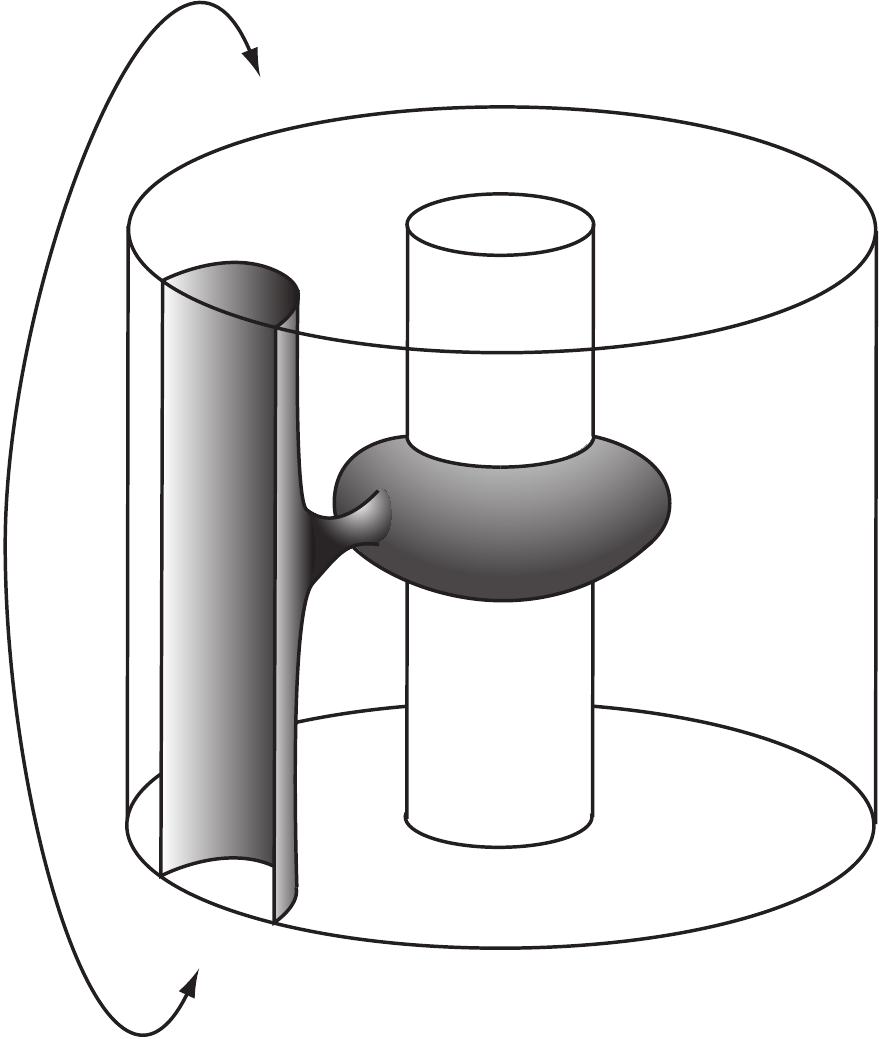}
        \caption{A toggle}
        \label{f:toggle}
        \end{center}
        \end{figure}

\end{proof}

\section{Toggles in Seifert fibered spaces}

The results of the previous section leave open the possibility that if $H$ can be isotoped to meet $T$ in a toggle then $H$ may be isotoped to meet $T$ in an arbitrarily large number of slopes. In the appendix R Wedimann shows that this can, indeed, happen. Here we prove that for ``most" Seifert fibered spaces it does not. In particular, we prove the following:

\begin{thm}
\label{t:SFSisotopy}
Let $M$ be a closed, totally orientable Seifert fibered space which is not a circle bundle with Euler number $\pm 1$. Let $H$ be a strongly irreducible Heegaard surface in $M$ and $T$ be a non-separating, vertical, essential torus. Then the isotopy class of $H$ determines at most two slopes on $T$.
\end{thm}

\begin{rmk}
The hypotheses of \fullref{t:SFSisotopy} can be relaxed to include any 3--manifold constructed in the following way. Begin with a Seifert fibered space which is not a circle bundle, with exactly two boundary components $T_1$ and $T_2$. Let $f_i$ denote a regular fiber on $T_i$. Construct $M$ by gluing $T_1$ to $T_2$ so that $|f_1 \cap f_2| \ne 1$. Let $T$ denote the image of $T_1$ and $T_2$ in $M$. If $H$ is any strongly irreducible Heegaard surface in $M$ then the conclusion of \fullref{t:SFSisotopy} holds for the pair $(T, H)$. 
\end{rmk}

\begin{proof}
Let $H_t$ be a compression free isotopy in which there are values $t_0$, $t_1$ and $t_2$ such $H_0$, $H_1$ and $H_2$ meet $T$ in different slopes (where $H_i=H_{t_i}$). Assume that $t_0$, $t_1$, and $t_2$ are consecutive with respect to this property, in the sense that there is no value $t \in (t_0,t_2)$ such that $H_t$ meets $T$ in some fourth slope. 

Let $N(T)$ denote a fibered, closed neighborhoood of $T$. Let $T_{\mu}$ and $T_{\nu}$ denote the boundary tori of $N(T)$. By \fullref{t:main} we know there is some $t_x \in (t_0,t_1)$ such that $H_{t_x}$ meets $N(T)$ in a toggle. Let $H_x=H_{t_x}$ and $\mu _x$ and $\nu _x$ denote the slopes of $H_x \cap T_{\mu}$ and $H_x \cap T_{\nu}$, respectively. Similarly, there is a $t_y \in (t_1,t_2)$ such that $H_{t_y}$ meets $N(T)$ in a toggle. Let $H_y=H_{t_y}$ and $\mu_y$ and $\nu _y$ denote the slopes of $H_y \cap T_{\mu}$ and $H_y \cap T_{\nu}$.

Let $M(T)$ denote the closure of $M \setminus N(T)$. As $H_x \cap
N(T)$ contains compressions on both sides of $H_x$ (the disks $D_1$
and $D_2$ from the proof of \fullref{t:main}) it follows from
strong irreducibilty that $H_x \cap M(T)$ is incompressible in
$M(T)$. By \cite[Theorem VI.34]{jaco:80} we may thus assume that each
component of $H_x \cap M(T)$ is horizontal or vertical. Similarly, we
may assume that each component of $H_y \cap M(T)$ is horizontal or
vertical.

We now show that $H_x \cap M(T)$ must be horizontal. First, note that since $T$ is non-separating $M(T)$ is connected. It follows that if $H_x \cap M(T)$ is not connected and one component is vertical then every component is vertical. This is because a horizontal component will meet every fiber, and hence will meet the fibers contained in the vertical components. We conclude the entire surface $H_x \cap M(T)$ is either vertical or horizontal. If it is vertical  then $\mu _x$ and $\nu_x$ will be fibers, and hence will represent the same slope on $T$. This contradicts the fact that they are on opposite sides of a toggle. 

We conclude the surface $H_x \cap M(T)$ is horizontal. We now assert that it consists of precisely two components, each with a single boundary component on each component of $\partial M(T)$. Suppose not. Then the surface $H_x \cap M(T)$ is a connected, horizontal surface. The two loops of $H_x \cap T_{\mu}$ inherit, from $H_x \cap M(T)$, orientations that agree on $T_\mu$. (In a totally orientable Seifert fibered space we can consistently orient each fiber. This defines a normal vector at every point of a horizontal surface.) Inspection of \fullref{f:toggle} indicates that these two loops inherit, from $H_x \cap N(T)$, orientations that disagree. As $H_x$ is orientable we have thus obtained a contradiction. 

A symmetric argument shows that $H_y \cap M(T)$ is a horizontal surface, made up of two components, each with one boundary component on each component of $\partial M(T)$. Any two horizontal surfaces in a Seifert fibered space differ by Dehn twists in vertical annuli and tori. (This is because given a spine $\Sigma$ of the base orbifold $\Sigma \times S^1$ cuts $M(T)$ into solid tori. As a horizontal surface intersects each such solid torus in meridian disks the only ambiguity arises from gluing the solid tori back together along vertical tori and annuli.) A Dehn twist in a vertical torus, however, does not change the boundary slopes of the surface. Similarly, as $M(T)$ is totally orientable a Dehn twist in a vertical annulus that has both boundary components on the same component of $M(T)$ will not change boundary slopes. We conclude that the pair $(\mu _y, \nu _y)$ can be obtained from the pair $(\mu _x, \nu _x)$ by Dehn twisting in annuli that have each of their boundary loops on different components of $\partial M(T)$. In other words, $(\mu _y, \nu _y)$ can be obtained from $(\mu _x, \nu _x)$ by {\it simultaneous} Dehn twisting along fibers. It follows that if $\mu _x=\mu _y$ then $\nu _x=\nu _y$, which is not the case, since by assumption there are exactly three distinct slopes among $\mu_x$, $\mu _y$, $\nu _x$, and $\nu _y$. 

We conclude, then, that $\mu _x=\nu _y$ or $\mu _y=\nu _x$. Without loss of generality assume the former. Now note that $\mu_y$ meets $\nu _y$ in a point, as one is obtained from the other by passing a toggle across $T$. Finally, this implies $\mu _x$ meets $\mu _y$ in a point. But $\mu _y$ is obtained from $\mu _x$ by Dehn twisting along a fiber. This can only happen if $\mu _x$ (and $\mu _y$) meets each fiber once. We conclude that each component of the horizontal surface $H_x \cap M(T)$ meets a regular fiber once, and hence $M(T)$ has no exceptional fibers (see \fullref{s:SFS} above).  

Finally, note that $M$ can be recovered from $M(T)$ by identifying its boundary components. But this must be done in such a way so that $\mu_x$ and $\nu _x$ meet in a point. Since these loops are at the boundary of a horizontal surface in $M(T)$, it must be the case that the Euler number of $M$ is positive or negative one.  
\end{proof}

\appendix

\section{Irreducible Heegaard splittings of circle bundles are unique
(by R Weidmann)}

The goal of this appendix is to prove \fullref{t:CircleBundles}. We denote the orientable circle bundle over the orientable surface $S_g$ of genus $g\ge 1$ with Euler number $e$ by $M_{g,e}$. 

In \cite{mr:98} Y~Moriah and J~Schultens show that all irreducible Heegaard splittings of Seifert manifolds are isotopic to horizontal or vertical Heegaard splittings. Moreover in the case of manifolds of type $M_{g,e}$ they show (see \cite{mr:98}, Corollary~0.5) that all irreducible Heegaard splittings of $M_{g,e}$ are vertical and of genus $2g+1$ if $e\neq\pm 1$ and horizontal of genus $2g$ if $e=\pm 1$. They further show that in the case $e \neq \pm 1$ the vertical splitting is unique up to isotopy. To prove \fullref{t:CircleBundles} it therefore suffices to show that all genus $2g$ horizontal Heegaard splittings of $M_{g,e}$ with $e=\pm 1$ are isotopic. The algebraic analogue of this statement is:

\begin{thm}
\label{t:Nielsen}
Let $M$ be an orientable circle bundle over an orientable surface of genus $g \ge 1$ with Euler number equal to $\pm 1$. Then any two generating tuples for $\pi_1(M)$ of cardinality $2g$ are Nielsen equivalent.
\end{thm}

We prove this theorem first, as the proof motivates the proof of \fullref{t:CircleBundles}.

Let $G$ be a group and $\mathcal T=(g_1,\ldots ,g_n)$ and $\mathcal T'=(g'_1,\ldots ,g'_n)$ be two tuples of elements. Recall that $\mathcal T$ and $\mathcal T'$ are called {\em elementary equivalent} if one of the following holds.

\begin{enumerate}
\item There exists some $\sigma\in S_n$ such that $g_i'=g_{\sigma(i)}$ for $1\le i\le n$.
\item $g_i'=g_i^{-1}$ and $g_j'=g_j$ for $j\neq i$.
\item $g_i'=g_ig_j^{\varepsilon}$ for some $i\neq j$ and $\varepsilon\in\{-1,1\}$. Furtermore $g_k'=g_k$ for $k\neq i$.
\end{enumerate}

We further say that $\mathcal T$ and $\mathcal T'$ are {\em Nielsen equivalent} if there exists a sequence of tuples $\mathcal T=\mathcal T_0,\ldots ,\mathcal T_k=\mathcal T'$ such that $\mathcal T_{i-1}$ and $\mathcal T_i$ are elementary equivalent for $1\le i\le k$.

\begin{proof}[Proof of \fullref{t:Nielsen}]
Let $g\ge 1$ and $e=\pm 1$. Note that 
$$\pi_1(M_{g,e})=\langle a_1,\ldots ,a_{2g},f\,|\,[a_1,f],\ldots ,[a_{2g},f],[a_1,a_2]\cdot \ldots\cdot[a_{2g-1},a_{2g}]f^{e}\rangle.$$
Let further $$p:\pi_1(M_{g,e})\to \pi_1(S_g)=\langle \bar a_1,\ldots ,\bar a_{2g}\,|\,[\bar a_1,\bar a_2]\cdot \ldots\cdot[\bar a_{2g-1},\bar a_{2g}]\rangle$$ be the projection given by $a_i\mapsto \bar a_i$ and $f\mapsto 1$. Recall that $\ker p=\langle f\rangle$. 

Note that $(a_1,\ldots ,a_{2g})$ is a generating tuple of $\pi_1(M_{g,e})$. 
To prove \fullref{t:Nielsen} it suffices to show that any generating tuple $(y_1,\ldots ,y_{2g})$ of $\pi_1(M_{g,e})$ is Nielsen equivalent to $(a_1,\ldots ,a_{2g})$.

A theorem of Zieschang \cite{Zie} states that in $\pi_1(S_g)$ any
is Nielsen equivalent to
$(\bar a_1,\ldots ,\bar a_{2g})$. It follows that for any generating
tuple\break $(y_1,\ldots ,y_{2g})$ of $\pi_1(M_{g,e})$ the tuple
$(p(y_1),\ldots ,p(y_{2g}))$ is Nielsen equivalent to $(\bar
a_1,\ldots ,\bar a_{2g})$. Thus $(y_1,\ldots ,y_{2g})$ and
$(a_1f^{z_1},\ldots ,a_{2g}f^{z_{2g}})$ are Nielsen equivalent for
some $z_i\in\mathbb Z$ for $1\le i\le 2g$.

It clearly suffices to show that for any $i=1,\ldots ,2g$ and
$\eta\in\{-1,1\}$ there exists a sequence of Nielsen equivalences that
replaces the tuple $(a_1f^{z_1},\ldots ,a_{2g}f^{z_{2g}})$ with
$(a_1f^{z_1},\ldots,a_{i-1}^{z_{i-1}},a_if^{z_i+\eta},
a_{i+1}^{z_{i+1}},\ldots\ldots ,a_{2g}f^{z_{2g}})$, ie, that replaces
$a_if^{z_i}$ with\break $a_if^{z_i+\eta}$ and leaves all other elements
unchanged. Note first that there is a cyclic conjugate $r$ of the
relator $[a_1,a_2]\cdot \ldots\cdot[a_{2g-1},a_{2g}]f^{e}$ if $\eta
=-e$ and of its inverse if $\eta=e$ such that (after using the fact
that $f$ commutes with the $a_i$) $$r=f^{-\eta}a_i^{-1}w_1a_iw_2$$
where $w_1$ and $w_2$ are words in $a_1,\ldots ,a_{i-1},a_{i+1},\ldots
,a_{2g}$ such that any of the $a_j$ ($j\neq i$) occurs in $w_1$ and
$w_2$ once with exponent $+1$ and once with exponent $-1$. In
particular we have the identity $a_if^{\eta}=w_1a_iw_2$ in $G$.

With appropriate Nielsen moves (left and right multiplication with the elements $a_jf^{z_j}$) we can replace $a_if^{z_i}$ with $w_1a_iw_2f^{z_i}$. (Note that the $f^{\pm z_j}$ cancel out as every $a_jf^{z_j}$ occurs once with exponent $+1$ and once with exponent $-1$.) As $w_1a_iw_2f^{z_i}=a_if^{\eta}f^{z_i}=a_if^{z_i+\eta}$ this proves the claim as all these Nielsen moves have left the $a_jf^{z_j}$ with $j\neq i$ untouched.
\end{proof}

We illustrate the main step of the above proof with an example. Suppose that $g=2$ and $e=1$ and that we want to show that $$(a_1f^{z_1},a_2f^{z_2},a_3f^{z_3},a_4f^{z_4})$$ is Nielsen equivalent to $$(a_1f^{z_1},a_2f^{z_2+1},a_3f^{z_3},a_4f^{z_4}),$$ this is the case with $i=2$ and $\eta=1$.

Note that the inverse of the long relation from the presentation of $\pi_1(M_{2,1})$ is $$f^{-1}a_4a_3a_4^{-1}a_3^{-1}a_2a_1a_2^{-1}a_1^{-1},$$ a cylic conjugate is $$a_2^{-1}a_1^{-1}f^{-1}a_4a_3a_4^{-1}a_3^{-1}a_2a_1.$$ As $f$ commutes with all $a_i$ we have the relation $$r=f^{-1}a_2^{-1}w_1a_2w_2$$ with $w_1=a_1^{-1}a_4a_3a_4^{-1}a_3^{-1}$ and $w_2=a_1$. Clearly $a_1$, $a_3$ and $a_4$ all occur twice in $w_1$ and $w_2$, once with exponent $+1$ and once with exponent $-1$.

It follows that by applying six Nielsen moves (where each takes one of the elements $a_1f^{z_1}$, $a_3f^{z_3}$, $a_4f^{z_4}$ or their inverse and multiplies the second element in the tuple from the left or right we can replace $(a_1f^{z_1},a_2f^{z_2},a_3f^{z_3},a_4f^{z_4})$ by 
\begin{gather*}(a_1f^{z_1},(a_1f^{z_1})^{-1}(a_4f^{z_4})(a_3f^{z_3})(a_4f^{z_4})^{-1}(a_3f^{z_3})^{-1}\hspace{1in}\\
\hspace{2in}a_2f^{z_2}(a_1f^{z_1}),a_3f^{z_3},a_4f^{z_4})\end{gather*}
All the $f^{\pm z_j}$ with $j\neq 2$ cancel in the second element of
the tuple, it follows that this new tuple is nothing but
\begin{gather*}(a_1f^{z_1},a_1^{-1}a_4a_3a_4^{-1}a_3^{-1}a_2f^{z_2}a_1,a_3f^{z_3},a_4f^{z_4})\hspace{1.7in}\\
=(a_1f^{z_1},w_1a_2w_2f^{z_2},a_3f^{z_3},a_4f^{z_4})=(a_1f^{z_1},a_2f\cdot f^{z_2},a_3f^{z_3},a_4f^{z_4})\\\hspace{2.4in}
=(a_1f^{z_1},a_2f^{z_2+1},a_3f^{z_3},a_4f^{z_4})\end{gather*}

\begin{proof}[Proof of \fullref{t:CircleBundles}]
As mentioned above, to prove \fullref{t:CircleBundles} it suffices to show that any two genus $2g$ horizontal Heegaard splittings of $M=M_{g,e}$ with $e=\pm 1$ are isotopic. We illustrate our proof of this assertion in \fullref{f:Richard}, in the case where the base orbifold of $M$ is a torus. The higher genus case is more difficult to see. 

        \begin{figure}[p]\def\strut{\vrule width0pt height6pt
depth 3pt}
\labellist\small\hair 6pt
\pinlabel {\vbox{\llap{Curve that\strut}\break\llap{bounds disk\strut}\break
\llap{in solid torus\strut}}} [r] at 195 716
\pinlabel {\vbox{\rlap{Slide handle through\strut}\break
\rlap{back face of cube\strut}}} [l] at 148 274
\pinlabel {\vbox{\rlap{isotopy\strut}}} [l] at 489 292
\pinlabel {\vbox{\rlap{Heegaard\strut}\break\rlap{surface\strut}}} 
[l] at 503 723
\pinlabel* {\vbox{\rlap{Dehn twist in torus\strut}\break\rlap{(right face
of cube)\strut}}} 
[l] <0pt,15pt> at 259 486
\pinlabel* {\vbox{\rlap{Pass spine across\strut}\break
\rlap{solid torus\strut}}} 
[l] <0pt,15pt> at 262 119
\pinlabel {Spine of Heegaard splitting} 
[t] at 99 377
\endlabellist
        \begin{center}
        \includegraphics[width=5 in]{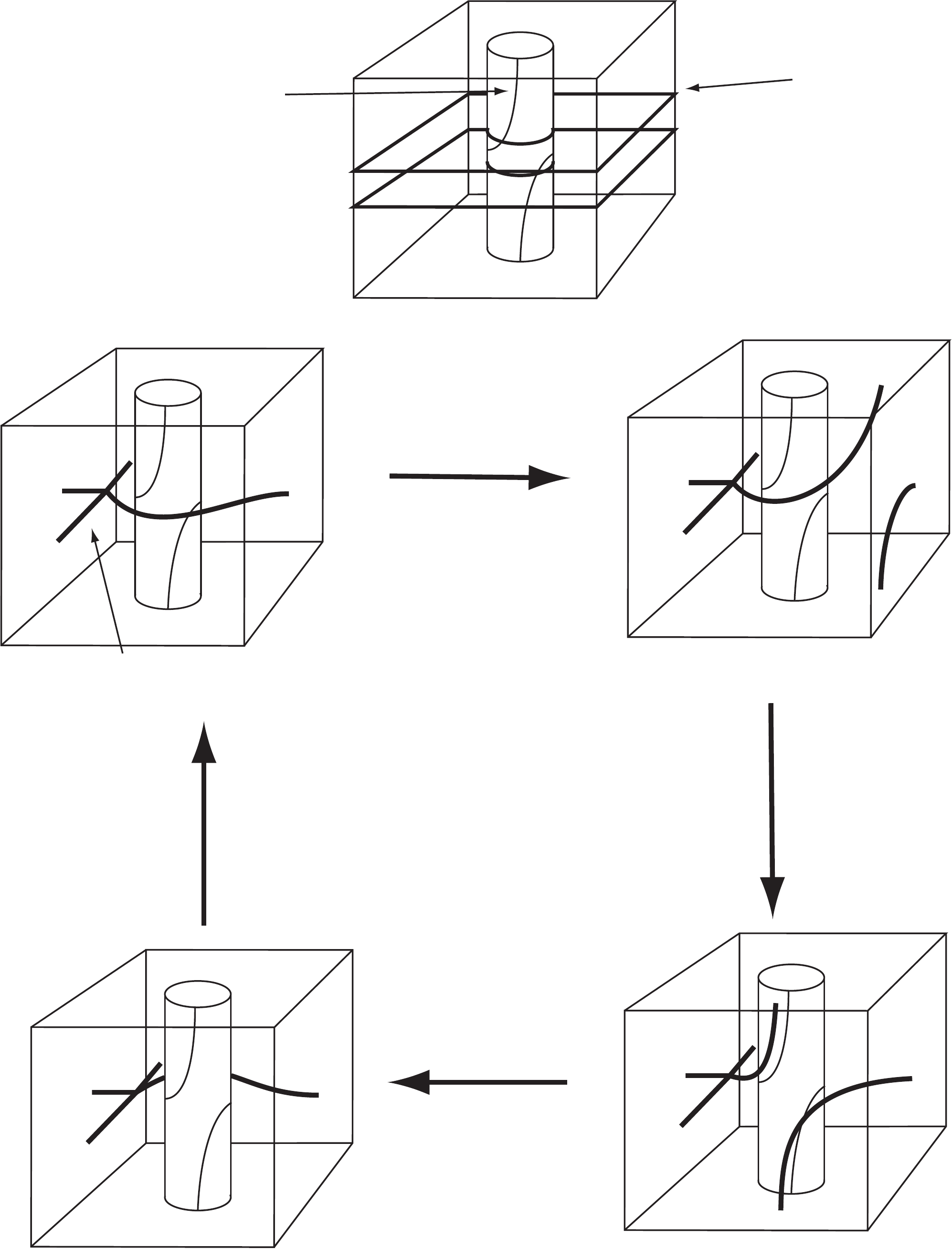}
        \caption{In $M_{1,\pm 1}$ Dehn twisting a horizontal Heegaard 
                 splitting about\break a vertical torus produces an isotopic 
                 splitting.}
        \label{f:Richard}
        \end{center}
        \end{figure}
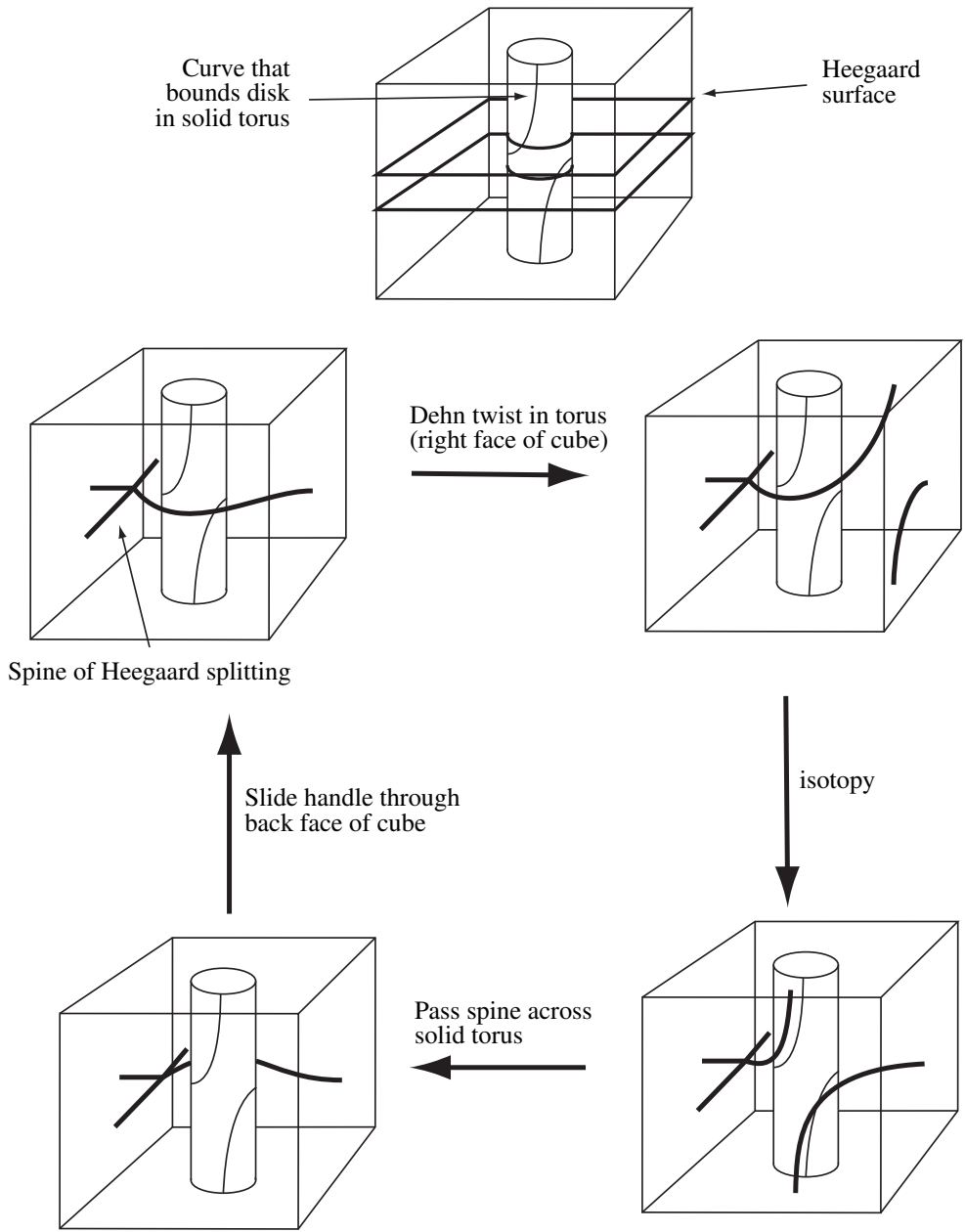

Let $f$ be a fiber of $M$. Note that $\hat M=M-N(f)\approx S_g^*\times S_1$ where $S_g^*\subset S_g$ is the once punctured orientable surface of genus $g$.  

Let further $\alpha_1,\ldots \alpha_{2g}\subset S_g^*$ be a canonical system of curves of $S_g$ with common base point $x$. Thus $\Gamma=\cup\alpha_i$ is a wedge of $2g$ circles and $S_g-\Gamma$ is a disk. Clearly we can assume that $\bar a_i=[\alpha_i]$ for $1\le i\le 2g$ if $$\pi_1(S_g,x)=\langle \bar a_1,\ldots ,\bar a_{2g}\mid[\bar a_1,\bar a_2],\ldots ,[\bar a_{2g-1},\bar a_{2g}]\rangle.$$
Let now $S$ be a horizontal Heegaard surface of genus $2g$. After an isotopy of $S$ we can assume that $S$ is horizontal at  the fibre $f$. Thus we can assume that $M=V\cup_{S}W$  where $V$ is the regular neighborhood in $\hat M$ of a horizontal surface. In particular there exists a Heegaard graph $\hat\Gamma$ (a core of $V$) that gets mapped homeomorphically to $\Gamma$ under the projection map $\pi:\hat M\to S_g^*\subset S_g$, $(x,z)\mapsto x$, in particular $\hat \Gamma$ is a wedge of $2g$ circles with single vertex $\hat x$ of valence greater than $2$. Denote the arc of $\hat\Gamma$ that gets mapped to $\alpha_i$ by $\beta_i$.

It is clear that the (ambient) isotopy class and in fact the homotopy  class of such a graph $\hat\Gamma$ is determined by the homotopy classes $[\beta_i]=a_if^{z_i}$ of the $\beta_i$ in $$\pi_1(M,\hat x)=\langle a_1,\ldots ,a_{2g},f\,|\,[a_1,f],\ldots ,[a_{2g},f],[a_1,a_2]\cdot \ldots\cdot[a_{2g-1},a_{2g}]f^{e}\rangle.$$ In particular we have a one to one correspondence between these types of Heegaard graphs and generating tuples $(a_1f^{z_i},\ldots, a_{2g}f^{z_{2g}})$ of $\pi_1(M,\hat x)$.

As in the proof of \fullref{t:Nielsen} it suffices to show that for any $i$ and $\eta\in\{-1,1\}$ the Heegaard graph $\hat\Gamma_1$ corresponding to the tuple $(a_1f^{z_i},\ldots, a_{2g}f^{z_{2g}})$ of $\pi_1(M,\hat x)$ is isotopic to the graph $\hat\Gamma_2$ corresponding to the tuple  $$(a_1f^{z_1},\ldots,a_{i-1}^{z_{i-1}},a_if^{z_i+\eta},a_{i+1}^{z_{i+1}},\ldots\ldots ,a_{2g}f^{z_{2g}}).$$

Choose a relation $r=f^{-\eta}a_i^{-1}w_1a_iw_2$ as in the proof of \fullref{t:Nielsen}. Recall that  $[\beta_j]=a_jf^{z_j}$ for $1\le j\le 2g$. Let $\hat w_1$ and $\hat w_2$ be the words obtained from $w_1$ and $w_2$ by replacing every occurence of  $a_j^{\pm 1}$ by $[\beta_j]^{\pm 1}$. As the fibre commutes with all $a_1$ it follows that we have the relation $$\hat r=f^{-\eta}[\beta_i]^{-1}\hat w_1[\beta_i]\hat w_2.$$  Let $\bar w_1$ and $\bar w_2$ be the path in $M_{g,e}$ obtained from $w_1$ and $w_2$ by replacing $a_i$ with $\beta_i$. Let further $\bar f$ be the fibre over $x$, clearly $[\bar f]=f$.

We have $[\bar w_i]=\hat w_i$ for $i=1,2$. Now there exists a map $h:D\to M_{g,e}$ of a disk $D$ such that $h(\partial D)$ is the path $\bar f^{-\eta}\beta_i^{-1}\bar w_1\beta_i\bar w_2$, such that $h$ is injective on the interior $D_0$ of $D$ and the projection onto the base space  maps $h(D_0)$ homeomorphically onto $S_g-\Gamma$. 

To see this note that any lift $\tilde\gamma$ of $\gamma=\bar f^{-\eta}\beta_i^{-1}\bar w_1\beta_i\bar w_2$ to the universal covering of $M$ is a simple closed curve that is contained in the boundary of the closure $\bar N$ of some component $\tilde N\approx D^2\times\mathbb R$ of the preimage of $N=M-\Gamma\times S_1\approx D^2\times S^1$ under the covering map. Clearly $\tilde\gamma$ bounds a properly embedded disk $D$ in $\bar N$. Note that $\tilde\gamma$ is transverse to the induced foliation of $\bar N$ by lines except in the subpath which is the lift of $\bar f^{-\eta}$. This subpath is contained in a leaf of the foliation. It follows that $D$ can be chosen to be horizontal in its interior, ie, transverse to the foliation. In particular the interior of $D$ intersects every line of the foliation of $\tilde N$ (the intrior of $\bar N$) exactly once. If follows that the restriction of the covering projection to $D$ is the desired map~$h$.

\begin{figure}[ht!]
\centerline{ \setlength{\unitlength}{.7cm}\small
\begin{picture}(15,6.5)
\put(7,3){$D$}
\put(5,5.5){\line(1,0){5}}
\put(7.5,5.5){\vector(1,0){.01}}\put(7.2,5.8){$\bar w_1$}
\put(5,5.5){\line(-1,-2){1.5}}\put(4.25,4){\vector(-1,-2){.01}}\put(3.7,4){$\beta_i$}
\put(10,5.5){\line(1,-2){1.5}}\put(10.75,4){\vector(1,-2){.01}}\put(10.9,3.9){$\beta_i$}
\put(11.5,2.5){\line(-3,-2){4}}\put(9.5,1.166){\vector(-3,-2){.01}}\put(9.6,.7){$\bar w_2$}
\put(3.5,2.5){\line(3,-2){4}}\put(5.5,1.166){\vector(3,-2){.01}}\put(4.6,.7){$\bar f^{\eta}$}
\bezier{6}(5,5.5)(5.05,5.3)(5.2,5.3)
\bezier{80}(9.8,5.3)(7.5,5.3)(5.2,5.3)
\bezier{70}(9.8,5.3)(10.55,3.9)(11.2,2.55)
\bezier{80}(7.55,0.1)(9.375,1.3)(11.2,2.55)
\bezier{6}(7.55,0.1)(7.5,.08)(7.5,-.16)
\end{picture}}
\end{figure}

We can now slide the edge $\beta_i$ over $\bar w_1$ and $\bar w_2$ which yields a new edge homotopic to $\beta_i\bar f^{\eta}$. No other edge is moved and the new edge can again be isotoped to map to $\alpha_i$. This proves the claim as the homotopy class of the new edge is  $[\beta_i\bar f^{\eta}]=[\beta_i][\bar f^{\eta}]=a_if^{z_i}f^{\eta}=a_if^{z_i+\eta}$.
\end{proof}

\bibliographystyle{gtart} \bibliography{link}

\end{document}